\documentclass[11pt,reqno,oneside]{amsart}

\usepackage[utf8]{inputenc}
\usepackage[T1]{fontenc}

\usepackage{amsmath,amsthm,amssymb}

\usepackage{graphics}
\usepackage[pagebackref]{hyperref}
\renewcommand*\backref[1]{\ifx#1\relax \else (Cited on page #1) \fi}

\usepackage{amsmath,amsthm,amssymb}

\usepackage{graphics}
\usepackage{hyperref}
\usepackage[dvipsnames]{xcolor}

\usepackage{numprint}
\npthousandsep{,}

\usepackage{booktabs}

\usepackage{mathtools}
\mathtoolsset{showonlyrefs}

\usepackage{enumitem}
\newlist{todolist}{itemize}{2}
\setlist[todolist]{label=$\square$}
\usepackage{pifont}
\usepackage[square,sort,comma,numbers]{natbib}

\definecolor{darkblue}{rgb}{0.0,0.0,0.3}
\hypersetup{colorlinks,breaklinks,
  linkcolor=darkblue,urlcolor=darkblue,
anchorcolor=darkblue,citecolor=darkblue}

\theoremstyle{plain}
\newtheorem*{theorem*}{Theorem}

\newtheorem*{proposition*}{Proposition}

\newtheorem*{corollary*}{Corollary}

\theoremstyle{definition}

\newtheorem*{remark*}{Remark}

\newtheorem*{algorithm*}{Algorithm}

\newtheorem*{experiment*}{Experiment}

\newcommand{\todaysdate}{\the\year.\ifnum\month<10 0\fi\the\month.\ifnum\day<10 0\fi\the\day}

\title[Studying number theory with deep learning]{Studying number theory with
deep learning:\\a case study with the Möbius and squarefree indicator functions}
\author{David Lowry-Duda}
\date{}

\begin{document}

\begin{abstract}
  Building on work of Charton, we train small transformer models to calculate
  the M\"{o}bius function $\mu(n)$ and the squarefree indicator function $\mu^2(n)$.
  The models attain nontrivial predictive power.
  We apply a mixture of additional models and feature scoring to give a
  theoretical explanation.
\end{abstract}

\maketitle

\section{Introduction}\label{sec:intro}

Many practical, state-of-the-art algorithms have come from deep learning.
In the last decade, Transformer-based models have been particularly prevalent.
Recent work has attempted to apply these techniques to pure mathematics,
including formal logic~\cite{hahn2021teaching, fengshi2021transformerbased},
symbolic integration and regression~\cite{rashidbarket2024transformers}, and
automated theorem proving~\cite{stanislaspolu2020generative}.
In each of these, machine learning models learn to manipulate abstract
mathematical symbols and perform symbolic computation.

Efforts to apply deep learning to concrete numerical calculations, such as basic
arithmetic operations, have been less successful.
Even basic arithmetic operations like multiplication or taking the modulus
after division appear very difficult to learn~\cite{kaiser2016neural,
palamas2017investigating}.
One promising line of research is due to Charton, who has used a
series of small (up to $6$ layers) sequence-to-sequence transformers to study
polynomial roots, matrix problems,
greatest common divisors, and more~\cite{charton2022polynomial,
charton2022linear, charton2024learning}.

In this paper, we investigate whether small transformers can learn functions
that are famously difficult to compute: the M\"{o}bius function $\mu(n)$ and the
squarefree indicator function $\mu^2(n)$.
The M\"{o}bius function assigns values based on the prime factorization of $n$ and
is defined by
\begin{equation}
  \mu(n) = \begin{cases}
  1 & \text{if } n=1, \\
  0 & \text{if } n \text{ has a squared prime factor}, \\
  (-1)^k & \text{if } n = p_1 p_2 \dots p_k \text{ where } p_i \text{ are distinct primes}.
  \end{cases}
\end{equation}
The squarefree indicator function, which is $0$ if $n$ is divisible by a
nontrivial square (we call such $n$ \emph{squarefull}) and is
otherwise $1$, is thus $\mu^2(n)$.

We apply the transformer architecture from Charton~\cite{charton2024learning} to
study $\mu(n)$ and $\mu^2(n)$.
Given Charton's successes in studying number-theoretic functions like GCDs, it
is natural to ask what transformers can learn about intricate multiplicative
functions.

\subsection{Setting expectations}\label{ssec:expect}

The obvious algorithm to compute $\mu(n)$ and $\mu^2(n)$ is to first factor $n$.
No known algorithm performs significantly better.
Adleman and McCurley~\cite{adleman1994openproblems} note that it is
unknown if there exists a polynomial-time algorithm to compute $\mu^2(n)$.
They also describe relationships between computing $\mu^2(n)$ and other
computationally hard problems in number theory.

In~\cite{shallit1985number}, Shallit and Shamir show that $\mu(n)$ can be
polynomially computed with an oracle that returns the number of divisors $d(m)$
for a single well-chosen integer $m$ (depending on $n$).
It is believed that computing the number of divisors is approximately as hard
as factoring, but this relationship is not proven.
Recently Booker, Hiary, and Keating gave a subexponential algorithm to compute
$\mu^2(n)$ assuming the generalized Riemann
Hypothesis~\cite{booker2015squarefree}; the running time is conjecturally slower
than current factoring algorithms, but is independent of factorization.

Closely related, there is a well-known but vague conjecture called ``M\"{o}bius
Randomness'' that approximately says that $\mu(n)$ does not strongly correlate
with any function that can be computed in polynomial time.
(See~\cite{sarnak2011lectures} for an introduction, or~\cite{green2012mobius} for
a more complete description of the same phenomenon).
Does M\"{o}bius Randomness imply that ML models should struggle to predict $\mu(n)$?
Or is it possible that ML could act like an approximate, probabilistic oracle?
This remains unknown.

\medskip{}

The probability that a random integer\footnote{Here and later, we interpret this
to mean the limiting probability as $X \to \infty$ after choosing an integer up
to $X$ uniformly at random. This is sometimes called \emph{natural density}.} is
squarefree is $1/\zeta(2) = 6/\pi^2 \approx 0.6079$.
The trivial algorithm of always guessing $1$ correctly predicts $\mu^2$
with probability $\approx 0.6079$.
The two cases $\mu(n) = \pm 1$ are approximately equally likely ($\approx
0.3039$), so the trivial algorithm for $\mu(n)$ is to always guess $0$.
This will be correct with probability $(1 - 1/\zeta(2)) \approx 0.3920$.

But it is not hard to do better.

A standard neural network trained on input-output pairs $(n, \mu^2(n))$ will
slowly learn to recognize if $4$ or $9$ or $25$ (and so on) divides $n$.
In these cases, it will correctly output $0$.
Guessing $1$ on the remaining cases will correctly compute $\mu^2(n)$ on a
random integer $n$ with probability approximately $0.967$.
Ellenberg makes this observation in his talks about machine learning (see for
example~\cite{ellenberg2024ml}).
Paraphrasing Ellenberg, continued training would likely lead to an algorithm
with nearly $100\%$ accuracy and nearly $0\%$ understanding.

In this paper, we seek learned insights \textbf{beyond
divisibility by small squares}.
To accomplish this, we choose an encoding of $n$ that deliberately obfuscates
square divisibility patterns (see \S\ref{sec:crt}).
Ultimately, we show that

\begin{enumerate}
  \item Choosing a representation beyond simply inserting $n$ and outputting
  $\mu(n)$ leads to different insights than divisibility by small squares.

  \item The statistical behavior of these predictions can be explained.
  Further, analyzing additional experiments with modified representations of
  the input $n$ confirm these explanations.

  \item As with previous experiments, model success in predicting $\mu(n)$
  reduces again to detecting squarefree and squarefull numbers; otherwise the
  models guess randomly.
\end{enumerate}

In the rest of this paper, we describe how we set up experiments to study the
M\"{o}bius function, found models with nontrivial predictive power, and then
iterated to explain the math underlying the models' success.

\section*{Acknowledgements}

This project grew out of the fall 2024 program on \emph{Mathematics and Machine
Learning} at the Harvard Center of Mathematical Sciences and Applications.
The ML component was performed on a system administered by Edgar Costa and Drew
Sutherland, using an ML architecture based on code generously shared by
Fran\c{c}ois Charton.
I also want to thank Angelica Babei, Edgar Costa, Mike Douglas, Jordan
Ellenberg, Noam Elkies, Xiaoyu Huang, Kyu-Hwan Lee, and Drew Sutherland for
their helpful discussions.
Finally, I also thank the Referees for their valuable remarks.

\medskip{}

This work was supported by the Simons Collaboration in Arithmetic Geometry,
Number Theory, and Computation via the Simons Foundation grant 546235.

\section{CRT Experiments}\label{sec:crt}

\defcitealias{int2int}{Int2Int}

We treat this as a supervised translation problem.
The setup is very similar to the setup of Charton in~\cite{charton2024learning},
and in particular uses code~\citepalias{int2int} that Charton wrote for experiments at
the fall 2024 \emph{Mathematics and Machine Learning} program at Harvard CMSA.\@
Code that replicates the primary experiments is available in an associated
code repository~\cite{lowryduda2025mobiusgithub}.
Integers $n$ were sampled uniformly randomly between $2$ and $10^{13}$,
represented as a sequence of tokens in ways that we make precise just below,
and used to train integer sequence-to-integer sequence transformers that learn
to output tokens representing either $\mu(n)$ or $\mu^2(n)$ by minimizing
cross-entropy between predictions and correct solutions.

Input integers $n$ are encoded as a sequence of $n \bmod p_j$ for the first
$100$ primes $p_j$,
\begin{equation}
  n \mapsto (n \bmod 2, n \bmod 3, \ldots, n \bmod 523, n \bmod 541).
\end{equation}
For each, the least nonnegative residue is chosen to represent the congruence
class (i.e.\ the common definition of \texttt{mod} in programming).
We refer to this as a Chinese Remainder Theorem (CRT) representation of $n$ using
the first $100$ primes.
Observe that any integer up to $\prod_{j = 1}^{100} p_j \lesssim 10^{219}$ can
be uniquely represented in such a way, more than sufficient to handle the
integers up to $10^{13}$ used for training and evaluation.

Each residue in the CRT representation is represented by the pair $(n \bmod
p_j, p_j)$, where the two integers are represented as sequences of digits in
base $1000$ using the sign `+' as a separator; as the largest modulus is $541$,
the vocabulary consisted of numbers between $0$ and $541$.
For example, the number $25$ would be encoded as the length $200$ sequence
\begin{equation}
  25
  \mapsto
  \text{`}
     \underbrace{\text{+ 1} \text{ + 2}}_{\bmod 2}
  \; \underbrace{\text{+ 1} \text{ + 3}}_{\bmod 3}
  \; \underbrace{\text{+ 0} \text{ + 5}}_{\bmod 5}
  \; \cdots
  \; \underbrace{\text{+ 25} \text{ + 523}}_{\bmod 523}
  \; \underbrace{\text{+ 25} \text{ + 541}}_{\bmod 541}
  \text{'}.
\end{equation}
Outputs were encoded as one of `- $1$', `+ $0$', or `+ $1$'.
During training, we experimented with other bases including base $4$, $10$, and
$25$; unlike~\cite{charton2024learning}, the choice of basis did not
have a strong effect.

The idea behind choosing the CRT representation is to deliberately make it
nontrivial to recognize divisibility by squares.
The only apparent way to determine if a CRT representation of $n$ is divisible
by $4$, say, is to completely reconstruct $n$ through the Chinese Remainder
Theorem, a nontrivial problem.

We trained transformers with $4$ layers, $256$ dimensions, and $8$ attention
heads using Adam~\cite{kingma2014adam} with inverse-square-root scheduling with
an initial learning rate of $2 \cdot 10^{-5}$ on batches of $64$ examples.
In total, our basic transformer model had just over \numprint{10000000}
parameters.
A total of $2 \cdot 10^6$ integers from $[2, 10^{13}]$ were uniformly sampled
with no repetition, and split randomly into a training set of size $1.8 \cdot
10^6$ and an evaluation set of size $2 \cdot 10^5$.
Training and testing sets were disjoint and did not change.
After each epoch ($\numprint{100000}$ examples), the models are tested on the
evaluation set.

\subsection*{Results of CRT experiments}

\begin{figure}[ht!]
  \parbox{0.49\textwidth}{
    \includegraphics[width=0.49\textwidth]{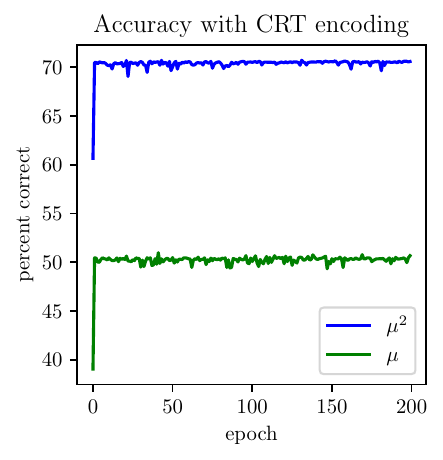}
  }\hfill
  \parbox{0.49\textwidth}{
    \includegraphics[width=0.49\textwidth]{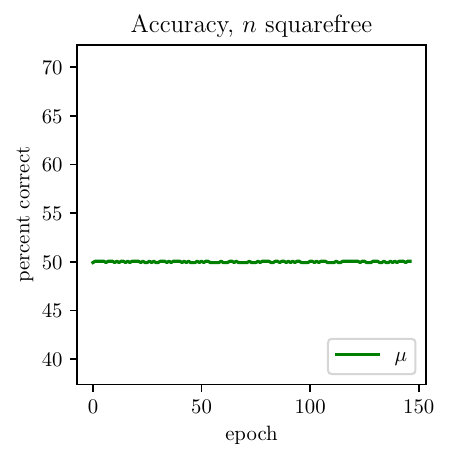}
  }
  \caption{Accuracy in predicting $\mu(n)$ and $\mu^2(n)$ from CRT
  representation. The left is for all $n$. The right is when restricted to
  squarefree $n$.}\label{fig:basic}
\end{figure}

Models very quickly learned to outperform trivial strategies.
Within $2$ epochs, models trained on CRT representations of
$\bigl(n, \mu^2(n)\bigr)$ correctly predict $70.64\%$ of examples in the
evaluation set (which should be compared to the default strategy, which is
correct $60.79\%$ of the time, cf.\ \S\ref{ssec:expect}).
Models trained on CRT representations of $\bigl(n, \mu(n)\bigr)$ have similar
behavior and achieve $50.98\%$ accuracy within the first $2$ epochs.
In the left plot in Figure~\ref{fig:basic}, we show the results for a typical
model --- small variations in the model parameters and setup all led to similar
behavior.
Variations were trained for many more epochs without improvement.

The prediction table in Table~\ref{tab:basic} hints at underlying behavior.
The best performing models for both $\mu(n)$ and $\mu^2(n)$ unsurprisingly tend
to overpredict the largest class (i.e.\ that $\mu(n) = 0$ and $\mu^2(n) = 1$),
leading to high accuracies for those outputs.
For $\mu(n)$, the model had very different behaviors at epochs with similar overall performance
concerning distinguishing between $\mu(n) = \pm 1$.
The models seem to pick randomly among these cases with no consistent behavior
throughout the training.

\begin{table}[htbp]
  \centering
  \begin{tabular}{cccccc}\toprule
  & \multicolumn{2}{c}{$\mu(n)$} & \multicolumn{2}{c}{$\mu^2(n)$} \\
  \cmidrule(lr){2-3}\cmidrule(lr){4-5}
  Output & \# in Eval Set & \# Recognized & \# in Eval Set & \# Recognized \\
  \midrule
  $-1$           & \numprint{6056} & \numprint{1203} &- &-\\
  $\phantom{-}0$ & \numprint{7852} & \numprint{6688} & \numprint{7836} & \numprint{4654} \\
  $\phantom{-}1$ & \numprint{6092} & \numprint{2305} & \numprint{12164} & \numprint{9473}\\
  \bottomrule\\
  \end{tabular}
  \caption{Number of correct predictions for $\mu(n)$ and $\mu^2(n)$.
  For each output, the total number of inputs $n$ yielding that output is
  recorded, as well as the number of those inputs correctly predicted by the
  model.
  }%
  \label{tab:basic}
\end{table}

As a follow-up experiment, we study CRT representations of
$\bigl(n, \mu(n)\bigr)$ with all training and testing data restricted to
squarefree integers.
(As $\mu^2(n)$ is trivial on this set, we do not consider this line of thinking
on $\mu^2$).
The generic behavior is shown on the right of Figure~\ref{fig:basic}.
Looking at this Figure, one might initially think there is no difference in
predictive power for $\mu$ when restricting to $n$ squarefree.
But the behavior is very different!
On the left, the model must choose between three candidate outputs, and the most
common choice is $\mu(n) = 0$ --- but this choice is removed as a possibility on
the right.
It is a numerical coincidence that the net accuracy is similar.
(See \S\ref{ssec:explanation} for a potential explanation).
This implies that all the predictive power in these models lies in their ability
to predict when $n$ is divisible by a square.

\section{Explaining Model Predictions}\label{sec:explain}

The models do not successfully distinguish between $\mu(n) = 1$ and $\mu(n) = -1$.
But how can we explain their success at identifying squarefree and squarefull
numbers?
In this section, study the models' behavior through two lines of analysis:
\begin{enumerate}
  \item We study how models perform on potential subproblems related to
  recovering $n$ from the Chinese Remainder Theorem description.
  \item We study the feature importance, and in particular study behavior
  before and after altering input datafiles.
\end{enumerate}

\subsection{False starts}\label{ssec:false}

One plausible explanation is that the models are learning some aspect of the
Chinese Remainder Theorem.
That is, given $(n \bmod p_j)$ for several primes $p_j$, the model successfully
recreates some representation of $n$.
We set up additional experiments, each using variations of the CRT
representation of $n$ as described in \S\ref{sec:crt} and with transformers of
similar architecture:

\begin{enumerate}
  \item
  \emph{Given $(n \bmod p_j)$ for the first $10$ primes $p_j$, output $n$.} ---
  The models completely fail.
  This is not a surprise, as this directly tests the Chinese Remainder Theorem
  and the output vocabulary is so large and non-repeating.

  \item
  \emph{Given $(n \bmod p_j)$ for the first $100$ primes $p_j$, and $n$
  restricted to $[1, 10^{13}]$, output if $n \in [1, 5 \cdot 10^{12}]$.
  Stated differently, learn the indicator function for an interval given a CRT
  representation.} ---
  The models perform no better than chance.

  \item
  \emph{Given $(n \bmod p_j)$ for the first $100$ primes, output $n \bmod 4$.
  Try to learn additional $2$-adic information.} ---
  The models perform poorly.
  We found this surprising as this is clearly related to determining whether $n
  \equiv 0 \bmod 4$ or not.

  \item
  \emph{Given $(n \bmod p_j)$ for the first $100$ primes except for
  $3$, try to learn $n \mod 3$.} ---
  The models perform no better than chance.

  \item
  \emph{Given $(n \bmod p_j)$ for the second hundred primes (from $547$
  to $1223$), try to compute $\mu^2(n)$.} ---
  The models perform no better than chance.
  We found this surprising and took this as a sign that the model is not
  using CRT reasoning.
\end{enumerate}

\subsection{Feature Importance}\label{ssec:limited_input}

\begin{figure}[ht!]
  \includegraphics[width=0.6\textwidth]{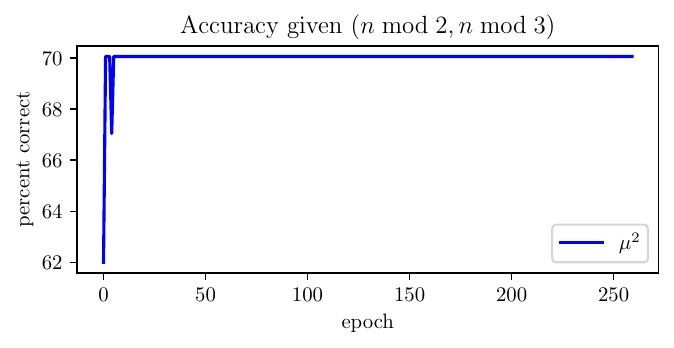}
  \caption{Accuracy in predicting $\mu^2(n)$ using only $(n \bmod 2, n \bmod 3)$
  as input.}\label{fig:limited}
\end{figure}

The failure of models to predict $\mu^2(n)$ given $(n \bmod p_j)$ for the
second $100$ primes $p_j$ leads one to naturally ask: \emph{which primes are the
most important?}
Using $(n \bmod 2, n \bmod 3)$ to predict $\mu^2(n)$ gives a model with the
behavior shown in Figure~\ref{fig:limited}.
This model predicts $\mu^2(n)$ with accuracy $70.1\%$, already significantly
better than chance and only slightly worse than the $70.64\%$ of the
best performing model from \S\ref{sec:crt}.

This gives a concrete, verifiable, purely mathematical claim: knowing $n \bmod
6$ is enough to guess $\mu^2(n)$ with probability about $70\%$.
This was not initially obvious, but we can now use this observation to explain
all model behavior.

\subsubsection*{Feature Corruption}

One method to determine whether the best-performing models really care only
about $(n \bmod 2, n \bmod 3)$ is through feature corruption:
take the models trained with the full CRT representation described in
\S\ref{sec:crt} and evaluate it on deliberately corrupted inputs.

We compare the performance across five different inputs:

\begin{enumerate}
  \item The \emph{true} input, in which no corruption has been introduced.
  \item Input in which $n \bmod 2$ is randomly altered to be $0$ or $1$ with
  equal probability.
  \item Input in which $n \bmod 3$ is randomly altered to be $0$, $1$, or $2$
  with equal probability.
  \item Input in which both $n \bmod 2$ and $n \bmod 3$ are randomly altered
  uniformly.
  \item Input in which $n \bmod 2$ and $n \bmod 3$ are preserved, but
  \emph{every other} $n \bmod p$ is randomly and uniformly altered.
\end{enumerate}

\begin{figure}[ht!]
  \parbox{0.49\textwidth}{
  {\small
    \begin{center}
      \begin{tabular}{ll}
        \toprule
        Data corruption experiment & Accuracy \\
        \midrule
        $\mu(n)$ true &  50.85 \\
        $\mu(n)$ with random $n \bmod 2$ & 38.49 \\
        $\mu(n)$ with random $n \bmod 3$ &  47.56 \\
        $\mu(n)$ with random $n \bmod 2, 3$ & 35.86 \\
        $\mu(n)$ only $n \bmod 2, 3$ correct &  50.90 \\
        \\
        $\mu^2(n)$ true & 69.72 \\
        $\mu^2(n)$ with random $n \bmod 2$ & 58.19 \\
        $\mu^2(n)$ with random $n \bmod 3$ &  67.77 \\
        $\mu^2(n)$ with random $n \bmod 2, 3$ &  55.17 \\
        $\mu^2(n)$, only $n \bmod 2, 3$ correct &  68.15 \\
        \bottomrule
      \end{tabular}
    \end{center}
  }
  }\hfill
  \parbox{0.49\textwidth}{
    \includegraphics[width=0.49\textwidth]{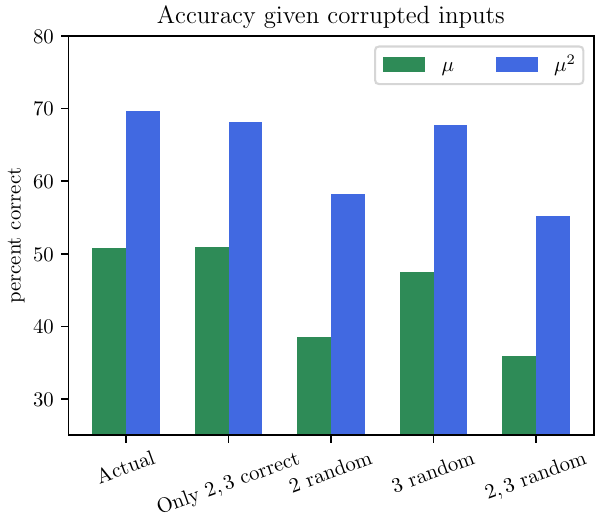}
  }
  \caption{Accuracy in predicting $\mu(n)$ and $\mu^2(n)$ from corrupted
  inputs. Random $n \bmod p$ means that in the CRT representation of $n$, the
  residue mod $p$ is replaced with a uniform random integer mod $p$.
  And ``only $n \bmod 2, 3$ correct'' means that the residues mod each prime $p
  \neq 2, 3$ were \emph{all} replaced with uniform random integers mod $p$. The
  plot on the right visualizes the data in the table on the left.
  }\label{fig:corrupt}
\end{figure}

For each of these inputs, the same set of $\numprint{20000}$ integers $n$ were
chosen and altered as described above.
The models were only evaluated (and not trained) on these inputs; hence no
information leaked between separate evaluations.
The result is shown in Figure~\ref{fig:corrupt}.

As is evident in Figure~\ref{fig:corrupt}, there is almost no change when
restricting to only use $n \bmod 2$ and $n \bmod 3$.
On the other hand, there is a large change when changing $n \bmod 2$, $n \bmod
3$, or both.
This shows that the model is most sensitive to $n \bmod 2$ and moderately
sensitive to $n \bmod 3$.
All other features are less important.

\subsection{Theoretical explanation}\label{ssec:explanation}

The model appears to be learning that numbers that are not squarefree are
probably divisible by a small square, and thus likely to be $0$ mod a small
prime. We explain this now.

To a subset of natural numbers $\mathcal{A}$, we associate a Dirichlet series
$D(s) = D_{\mathcal{A}}(s)$ as a generating function,
\begin{equation}
  D_\mathcal{A}(s) = \sum_{n \in \mathcal{A}} \frac{1}{n^s}.
\end{equation}
An application of Perron's formula (see e.g.~\cite{montgomery_multnumtheory})
shows that if $D(s)$ has a pole of order $\leq 1$ at $s = 1$, then the residue
at $s = 1$ gives the natural density of $\mathcal{A} \subset \mathbb{N}$.
We can construct the generating Dirichlet series for the relevant subsets of natural
numbers and compute their residues in terms of the Riemann zeta function
$\zeta(s) = \sum_{n \geq 1} 1/n^s$ below.

By unique factorization, the series associated to squarefree numbers is
\begin{equation}
  \prod_{p \text{ prime}}
  \Bigl( 1 + \frac{1}{p^s} \Bigr)
  =
  \frac{\zeta(s)}{\zeta(2s)}.
\end{equation}
The pole at $s = 1$ has residue $1/\zeta(2) = 6/\pi^2$, confirming
that the probability a random integer is squarefree is $6/\pi^2$.

The Dirichlet series associated to odd and even squarefree numbers are,
respectively,
\begin{equation}
  \prod_{\substack{p \text{ prime} \\ p \neq 2}}
  \Bigl( 1 + \frac{1}{p^s} \Bigr)
  =
  \frac{(1 - 1/2^s) \zeta(s)}{(1 - 1/4^s) \zeta(2s)},
  \quad
  \frac{1}{2^s}
  \prod_{\substack{p \text{ prime} \\ p \neq 2}}
  \Bigl( 1 + \frac{1}{p^s} \Bigr)
  =
  \frac{1}{2^s}
  \frac{(1 - 1/2^s) \zeta(s)}{(1 - 1/4^s) \zeta(2s)},
\end{equation}
from which it follows that the ratio of positive integers up to $X$ that are odd
(resp.\ even) and squarefree tends to $\frac{2}{3} \frac{6}{\pi^2}$
(resp.\ $\frac{1}{3} \frac{6}{\pi^2}$).
In particular, a squarefree integer is twice as likely to be odd as even.

For this classification problem, we're interested in the converse.
Writing $P(\textup{sqfree})$ to mean the probability that a random integer is
squarefree, and similarly for $P(\textup{even})$ and $P(\textup{odd})$, then
basic probability shows that
\begin{align}
  P(\textup{sqfree} | \textup{even})
  =
  \frac{P(\textup{even and squarefree})}{P(\textup{even})}
  =
  \frac{\frac{1}{3} \frac{6}{\pi^2}}{\frac{1}{2}}
  &\approx
  0.4052 \\
  P(\textup{sqfree} | \textup{odd})
  =
  \frac{P(\textup{odd and squarefree})}{P(\textup{odd})}
  =
  \frac{\frac{2}{3} \frac{6}{\pi^2}}{\frac{1}{2}}
  &\approx
  0.8105.
\end{align}

This already leads to a better-than-random strategy for both predicting $\mu(n)$
and $\mu^2(n)$: if $n$ is even, guess that $\mu(n) = 0$ (correct with
approximate probability $1 - 0.4052 \approx 0.5948$); otherwise, guess $\mu^2(n)
= 1$ (correct $\approx 0.8105$) and $\mu(n) = 1$ (correct $\approx 0.4052$).
In total, this would predict $\mu^2(n)$ correctly with approximate probability
$0.5 \cdot 0.5948 + 0.5 \cdot 0.8105 \approx 0.7026$, and would predict $\mu(n)$
correctly with approximate probability $0.5 \cdot 0.5948 + 0.5 \cdot 0.4052
\approx 0.5$.
These both improve upon the trivial algorithms from \S\ref{ssec:expect} and are
only moderately worse than the best CRT model in \S\ref{sec:crt}.

More generally, if $\{ p_1, \ldots, p_N\}$ and $\{q_1, \ldots, q_M\}$ are two
disjoint sets of primes, then the density of squarefree numbers that are
divisible by each of the $p_i$ and that are not divisible by any of the $q_j$ is
\begin{equation}
  \prod_{p_i} \Bigl( \frac{1}{p_i + 1} \Bigr)
  \prod_{q_j} \Bigl( \frac{q_j}{q_j + 1} \Bigr)
  \frac{6}{\pi^2},
\end{equation}
from which it follows that the probability that such an $n$ is squarefree given
these divisibility constraints is
\begin{equation}
  \prod_{p_i} \Bigl( \frac{1}{p_i + 1} \Bigr)
  \prod_{q_j} \Bigl( \frac{q_j^2}{q_j^2 - 1} \Bigr)
  \frac{6}{\pi^2}.
\end{equation}

Straightforward but tedious combinatorial analysis shows that using only
divisibility by the first $25$ primes leads to a strategy to approximate
$\mu^2(n)$ with accuracy $70.34\%$.
Combinatorial difficulties make it difficult to exactly compute the effect of
the cross correlations for the $100$ primes, but it is consistent with the
models' behavior.

The models are learning differences in these conditional probabilities from the
special cases when $n \equiv 0 \bmod p_j$ for one of the primes $p_j$.
This behavior is similar to the base-distinguishing behavior
in~\cite{charton2024learning}: there, models that predict GCDs would quickly
learn to determine divisibility rules for multiples of the base used in
transformer encoding by recognizing whether the rightmost digits in the base are
$0$.
The CRT representation acts as though we use several bases simultaneously.

\section{Discussion}\label{sec:discussion}

Previous neural network experiments trained to predict $\mu^2(n)$ from $n$
appear to learn the is-divisible-by-a-small-prime-square function.
Small transformer models trained on CRT representations of $n$ to predict
$\mu^2(n)$ appear to learn the is-divisible-by-a-small-prime function.
It is retrospectively clear that ML models should be able to learn
cross-correlations and make improved predictions based on non-uniformities.
But it was not obvious what the behavior here would be, or that this
analysis would work so well.

A fundamental difficulty came from a difference in perspective: the model is
trying to \textbf{predict} $\mu(n)$ and $\mu^2(n)$, not to \textbf{compute}
$\mu(n)$ and $\mu^2(n)$.
Discussing the behavior of these models with other mathematicians often led to
considering algorithms to \emph{compute} $\mu(n)$ given bits of information.
This was already conceptually misleading, obfuscating what the model was
actually doing.

Small transformers like Int2Int~\cite{int2int} are rapidly trainable, allow a
variety of input and output options, and are useful to act as \textbf{one-sided
oracles} to determine whether inputs contain enough information to evaluate an
output.
It's ``one-sided'' because the model might simply fail to approximate the
function, and this failure cannot distinguish between insufficient data or
insufficient model size or insufficient informational content.
The model acts like an ``oracle'' because there are no explanations for the
insights given.
Nonetheless, if the model makes successful predictions, this indicates that
there is nontrivial information and correlation between the inputs and outputs.

To end, we issue a challenge.
In the set of experiments and analysis described here, we ultimately had to
consider only $\mu^2(n)$ as the models were unable to distinguish between
$\mu(n) = 1$ and $\mu(n) = -1$.
Is it possible to do better?

\begin{quote}
\textbf{M\"{o}bius Challenge}\\
Train an ML model with inputs computable in time $\ll \log^A(n)$ for some
finite $A$ that distinguishes between $\mu(n) = 1$ and $\mu(n) = -1$ with
probability greater than $51\%$.
\end{quote}

\bibliographystyle{alpha}
\bibliography{bibfile}
\end{document}